\documentclass[11pt]{article}
\usepackage{bbm}
\usepackage{mathrsfs}
\usepackage{amsmath}
\usepackage{amsfonts}
\usepackage{graphicx}
\usepackage{amssymb}
\usepackage{cite}
\usepackage{cases}
\usepackage{empheq}
\usepackage{subeqnarray}
\usepackage[titletoc]{appendix}

\newtheorem{condition**}{A*}
\newtheorem{condition***}{C*}
\newtheorem{condition*}{C}

\newtheorem{proposition}{Proposition}[section]

\newtheorem{theorem}{Theorem}[section]
\newtheorem{lemma}{Lemma}[section]
\newtheorem{remark}{Remark}[section]

\newenvironment{keywords}{{\bf Key words: }}{}

\begin{document}

\title{A kind of linear quadratic non-zero sum differential game of backward stochastic differential equation with asymmetric information
}

\author{Guangchen Wang\footnotemark[1]\quad Hua Xiao\footnotemark[2] \quad Jie Xiong\footnotemark[3]
}
\renewcommand{\thefootnote}{\fnsymbol{footnote}}
\footnotetext[1]{School of Control Science and Engineering, Shandong
University, Jinan 250061, China (wguangchen@sdu.edu.cn). This author
acknowledges the support in part from the NSF of China under Grants
11371228, 61422305 and 61304130, by the NSF for Distinguished Young
Scholars of Shandong Province of China under Grant JQ201418, by the
Program for New Century Excellent Talents in University of China
under Grant NCET-12-0338, and by the Research Fund for the Taishan
Scholar Project of Shandong Province of China.}
\footnotetext[2]{Corresponding Author. School of Mathematics and
Statistics, Shandong University, Weihai 264209, China
(xiao\_\;hua@sdu.edu.cn). This author acknowledges the financial
support from the NSF of China under Grants 11471192 and 61573217.}
\footnotetext[3]{Department of Mathematics, University of Macau,
Taipa Macau, China (jiexiong@umac.mo, jxiong@math.utk.edu). This author acknowledges the financial support from  MYRG2014-00015-FST.}


\maketitle

\begin{abstract}
This paper focuses on a kind of linear quadratic non-zero sum differential game
driven by backward stochastic differential equation with asymmetric
information, which is a natural continuation of \cite{WY2010,WY2012}.
Different from \cite{WY2010,WY2012}, novel motivations for studying this kind of game
are provided. Some feedback Nash equilibrium points are uniquely obtained by
forward-backward stochastic differential equations, their filters and the corresponding Riccati equations with Markovian setting.
\end{abstract}



\section{Introduction}

Stochastic differential game plays an important role in lots of fields.
Many researchers investigated this problem under various setups
\cite{BSYY2014, ES2011, OS2014, Yong2002}. Recently,
\cite{WY2010} studied a non-zero sum differential game
of nonlinear backward stochastic differential equation (BSDE, for short).
Later, in \cite{WY2012}, they generalized the game in
\cite{WY2010} to the partial information case, and obtained an
open-loop Nash equilibrium point for a linear quadratic (LQ, for short) game
with same observable information. In some situations of real
markets, say, insider trading, one investor may get more
information than the others, and then, this investor can make a
better decision than the others. It implies that asymmetric
information has effect on the decision making. Such a kind of effect
is pervasive in reality, but is usually ignored in literature. To
fill in the gap, this paper initiates the study of an LQ non-zero
sum differential game of BSDE with asymmetric information. This study
can be regarded as a first step to investigate such a
kind of differential game with asymmetric information.

This paper is closely related to \cite{CX2014,SWX2015}, where the state satisfies a (forward) stochastic differential equation (SDE, for short), and thus the BSDE appears as an adjoint of the state equation. In this paper, the state is governed by a BSDE rather than an SDE. Since the construction and property of BSDE are
essentially different from those of SDE, the game of BSDE captures different scenarios. See, e.g., Section 2.1 for more information. This paper is also related to \cite{Hamadene1999,
LZ2001, MY2006, Yu2012,YJ2008,Zhang2011, HX2014}, where asymmetric
information is not considered. Therefore, this paper is
distinguished from the exiting references about stochastic
differential game.

The rest of this paper is organized as follows. In Section \ref{Section1.2}, a kind of LQ game
of BSDE with asymmetric information is formulated in detail and an
open-loop Nash equilibrium point is derived. Section \ref{SC} is
devoted to solving three concrete cases of the LQ game. Feedback Nash
equilibrium points are uniquely obtained by the filters of
forward-backward SDEs (FBSDEs, for short). One numerical example is also shown. In Section \ref{CR}, some
concluding remarks are given. Finally, in Appendix, several examples are shown to illustrate that the
special cases we studied in Section \ref{SC} are realistic.



 \section{Problem formulation and equilibrium points}\label{Section1.2}

Let us begin with a complete filtered probability space $(\Omega,
\mathscr{F}, (\mathscr{F}_{t})_{0\leq t\leq T}, \mathbb{P})$, in
which $\mathscr{F}_{t}$ denotes a natural filtration generated by a
two dimensional standard Brownian motion $w(t)=(w_{1}(t),
w_{2}(t))^{*}$. Suppose that $\mathscr{F}=\mathscr{F}_{T}$, $\mathbb{E}$
is the expectation with respect to $\mathbb{P}$, and $T>0$ is a
fixed time horizon. We denote by the superscript $*$ the transpose
of vectors or matrices, by $|\cdot|$ the norm, and by
$\mathscr{F}_{t}^{X}$ the filtration generated by a stochastic
process $X$, i.e., $\mathscr{F}_{t}^{X}=\sigma\big\{X(s), 0\leq
s\leq t\big\}.$ We call
$\mathbb{E}\left(h(t)|\mathscr{F}_{t}^{X}\right)$ the optimal filter
of $h(t)$ with respect to $\mathscr{F}_{t}^{X}$. We also give the notations $\tilde{h}(t)=\mathbb{E}\left(h(t)|\mathscr{F}_{t}^{w_{2}}\right)$
and $\hat{h}(t)=\mathbb{E}\left(h(t)|\mathscr{F}_{t}^{w_{1}}\right)$.

Let $\mathscr{G}^{i}_{t}\subseteq\mathscr{F}_{t}$ be a given
sub-filtration, which represents the information available to the
player $i\, (i = 1, 2)$ up to the time $t$. If
$\mathscr{G}^{i}_{t}=\mathscr{F}_{t}$ (resp.
$\mathscr{G}^{i}_{t}\subset \mathscr{F}_{t}$), we call the
information available to the player $i$ \emph{complete} (resp.
\emph{partial}). If $\mathscr{G}^{1}_{t}\neq\mathscr{G}^{2}_{t}$
(resp. $\mathscr{G}^{1}_{t}=\mathscr{G}^{2}_{t}$ ), we call the
information available to two players \emph{asymmetric} (resp.
\emph{symmetric}). For simplicity, we usually omit the terminology
``complete information''.

\subsection{An economic example}\label{EE}

Suppose that a consumer has a reward $\xi>0$ at the terminal time $T$ and continuously consumes between 0 and $T$. Here $\xi$ is an
$\mathscr{F}_{T}$-measurable and square-integrable random variable. Let $c_{1}(t)$ and $c_{2}(t)$ be the consumption rates about two different consumables $F_{1}$ (such as certain kind of meat) and $F_{2}$ (such as certain kind of vegetable), respectively. Let $p_{1}(t)$ and $p_{2}(t)$ be the prices of  $F_{1}$ and $F_{2}$, respectively, which are $\mathscr{F}_{t}$-adapted processes.
Set $\mathscr{P}^{i}_{t}=\sigma\{p_{i}(s); 0\leq s\leq t\}$, and let
\begin{align*}
  &\mathscr{C}_{i}=\left\{c_{i}(\cdot)|\, c_{i}(t) \mbox{ is } \mathscr{G}^{i}_{t}  \hbox{-adapted and square-integrable} \right\}
\end{align*}
be the set of all consumption rates $c_i(t)$, where $\mathscr G^i_t\subseteq\mathscr P^i_t$ $(i = 1, 2)$. It implies that the consumer chooses $c_{i}(t)$ depending on $\mathscr{G}^{i}_{t}$ $(i = 1, 2)$. This is reasonable in reality.

Let $y^{c_{1}, c_{2}}(t)$ be the Kreps-Porteus recursive utility of the consumer.
According to \cite{Peng1997}, a special case of $y^{c_{1}, c_{2}}(t)$ is modeled by
$$
\left\{
\begin{aligned}
-dy^{c_{1}, c_{2}}(t)&\ = \big(c_{1}(t)+c_{2}(t)-y^{c_{1}, c_{2}}(t)\big)dt-z_{1}^{c_{1}, c_{2}}(t)dw_{1}(t)-z_{2}^{c_{1}, c_{2}}(t)dw_{2}(t), \\
y^{c_{1}, c_{2}}(T)&\ =\xi.
\end{aligned}
\right.
$$
Define the performance functional as
$$\begin{aligned}\label{UF}
 J_{i}\big(c_{1}(\cdot), c_{2}(\cdot)\big)=\mathbb{E}\left[\int_{0}^{T}\frac{1}{2} \left(c_{i}(t)-e_{i}(t)\right)^{2}dt-r_{i} y^{c_{1}, c_{2}}(0)\right],
\end{aligned}
$$
where $e_{i}$ is a deterministic and uniformly bounded function, and is interpreted as a dynamic benchmark; $r_{i}$ is a positive constant $(i=1, 2)$. It is natural that the consumer wants not only to prevent $c_{i}(t)$ from large deviation, but also to maximize $y^{c_{1}, c_{2}}(0)$. That is,
$$
 \left\{
 \begin{aligned}
& J_{1}(c_{1}^{*}(\cdot), c_{2}^{*}(\cdot)) = \min\limits_{c_{1}(\cdot)\in \mathscr{C}_{1}} J_{1}(c_{1}(\cdot), c_{2}^{*}(\cdot)),\\
& J_{2}(c_{1}^{*}(\cdot), c_{2}^{*}(\cdot)) = \min\limits_{c_{2}(\cdot)\in \mathscr{C}_{2}} J_{2}(c_{1}^{*}(\cdot), c_{2}(\cdot)).
\end{aligned}
 \right.
$$
Note that $y^{c_1, c_2}(t)$ satisfies a BSDE and $\mathscr G^1_t$ is not always equal to $\mathscr G^2_t$.
Then the economic example can be regarded as a special LQ non-zero sum differential game of BSDE with asymmetric information.

\subsection{Problem formulation}\label{PF}            

Motivated by the above example, we consider the controlled linear BSDE
\begin{equation}\label{Eq1}
\left\{\begin{aligned}
-dy^{v_{1}, v_{2}}(t) = & \Big(a(t)y^{v_{1}, v_{2}}(t)+b_{1}(t)v_{1}(t)+b_{2}(t)v_{2}(t)+\sum_{j=1}^{2}f_{j}(t)z_{j}^{v_{1}, v_{2}}(t)+c(t)\Big)dt\\
                            &\ -z_{1}^{v_{1}, v_{2}}(t)dw_{1}(t)-z_{2}^{v_{1}, v_{2}}(t)dw_{2}(t), \\
y^{v_{1}, v_{2}}(T) = & \xi,
\end{aligned}
\right.
\end{equation}
and the cost functional
\begin{equation}\label{Eq2}
\begin{aligned}
\mathcal{J}_{i}\big(v_{1}(\cdot), v_{2}(\cdot)\big)=\frac{1}{2} \mathbb{E}\left\{ \int_{0}^{T}\left[l_{i}(t)\big(y^{v_{1}, v_{2}}(t)-k_{i}(t)\big)^{2}+m_{i}(t)\big(v_{i}(t)-n_{i}(t)\big)^{2}\right]dt\right.\\
\left.+r_{i}\big(y^{v_{1}, v_{2}}(0)-h_{i}\big)^{2}\vphantom{\int_{0}^{T}}\right\} \quad (i=1, 2).
\end{aligned}
\end{equation}
 Here $a, b_{1}, b_{2}, f_{1}, f_{2}, c, k_{1}, k_{2}, n_{1}$ and
$n_{2}$ are uniformly bounded and $\{\mathscr{F}_{t}, 0\leq t\leq
T\}$-adapted; $h_{1}$ and $h_{2}$ are given constants;  $l_{1}$,
$l_{2}$, $ m_{1}$ and $m_{2}$ are positive, uniformly bounded and
$\{\mathscr{F}_{t}, 0\leq t\leq T\}$-adapted; $r_{1}$ and $r_{2}$
are two nonnegative constants; $\xi$ is an
$\mathscr{F}_{T}$-measurable and square-integrable random variable;
$v_{1}(\cdot)$ and $v_{2}(\cdot)$ are the control processes of the
player 1 and the player 2, respectively. We use the
notation $(y^{v_{1}, v_{2}}, z_{1}^{v_{1}, v_{2}}, z_{2}^{v_{1},
v_{2}})$ to denote the dependence of the state on the control
$(v_{1}, v_{2})$. Introduce the admissible control set for the
player $i \ (i = 1, 2)$
\begin{align*}
  &\mathscr{U}_{i}=\left\{v_{i}(\cdot)|\, v_{i}(t) \mbox{ is } \mathscr{G}^{i}_{t}  \hbox{-adapted and square-integrable} \right\}\,.
\end{align*}
Each element of $\mathscr{U}_{i}$ is called an open-loop
admissible control for the player $i\, (i = 1, 2)$. $\mathscr{U}_{1}\times\mathscr{U}_{2}$ is the set of open-loop admissible controls for the players. Suppose that each player $i$ hopes to minimize her/his cost functional
$\mathcal {J}_{i}(v_{1}(\cdot), v_{2}(\cdot))$ by selecting a suitable admissible control
$v_{i}(\cdot)\, (i = 1, 2)$. Then the problem is to look for
$(u_{1}(\cdot), u_{2}(\cdot)) \in \mathscr{U}_{1}\times\mathscr{U}_{2}$, which is called a Nash equilibrium point of the game, such that
$$
\label{Eq39}
\left\{\begin{aligned}
& \mathcal{J}_{1}(u_{1}(\cdot), u_{2}(\cdot)) = \min\limits_{v_{1}(\cdot)\in \mathscr{U}_{1}} \mathcal{J}_{1}(v_{1}(\cdot), u_{2}(\cdot)),\\
& \mathcal{J}_{2}(u_{1}(\cdot), u_{2}(\cdot)) = \min\limits_{v_{2}(\cdot)\in \mathscr{U}_{2}} \mathcal{J}_{2}(u_{1}(\cdot), v_{2}(\cdot)),
\end{aligned}
\right.
$$
subject to (\ref{Eq1}) and (\ref{Eq2}). We call the game problem an LQ non-zero sum stochastic differential
game of BSDE with asymmetric information. For simplicity, we denote
the problem by \textbf{Problem (AI)}, and abbreviate $(y^{u_{1},
u_{2}}, z_{1}^{u_{1}, u_{2}}, z_{2}^{u_{1}, u_{2}})$ by $(y, z_{1}, z_{2})$. Clearly, Problem (AI) covers the example in Section \ref{EE} as a special case.

The main goal of this paper is to derive some Nash equilibrium points in the feedback form of the filtered states. However, since $\mathscr G^i_t$ available to the player $i\, (i=1, 2)$ is only an abstract sub-filtration of $\mathscr F_t$, it is impossible to obtain feedback Nash equilibrium points in general. Then some special information structures for $\mathscr G^i_t\, (i=1, 2)$ are desirable to
reach the goal. For example,
(i) $\mathscr{G}^{1}_{t}=\mathscr{G}^{2}_{t}=\mathscr{F}_{t}^{w_{2}},$ i.e., two players have access to the same observation information;
(ii) $\mathscr{G}^{1}_{t}=\mathscr{F}_{t}$ and $\mathscr{G}^{2}_{t}=\mathscr{F}^{w_{2}}_{t},$ i.e., one player has more information at any time than the other player;
(iii) $\mathscr{G}^{1}_{t}=\mathscr{F}^{w_{1}}_{t}$ and $\mathscr{G}^{2}_{t}=\mathscr{F}^{w_{2}}_{t},$ i.e., two players have independent
observation information and do not share all of their information with each other.
These special information structures are inspired by Remarks A.1, A.2 and A.3 in Appendix, respectively.

\subsection{Nash equilibrium point}\label{EP}    

The following proposition is an immediate result of Theorem 2.1 in \cite{WY2012}. It is very helpful for us to discuss some details and special cases of Problem (AI).

\begin{proposition}\label{Proposition3.3}
$(u_{1}, u_{2})$ is a Nash equilibrium point of Problem (AI) if and only if
$(u_{1}, u_{2})$ is in the form of
\begin{equation}\label{Eq57}
 \left\{
 \begin{aligned}
 u_{1}(t)=\frac{\mathbb{E}\left(b_1(t)x_{1}(t)|\mathscr{G}^{1}_{t}\right)}{\mathbb E\left(m_1(t)|\mathscr{G}^{1}_{t}\right)}+\frac{\mathbb{E}\left(m_1(t)n_{1}(t)|\mathscr{G}^{1}_{t}\right)}{\mathbb E\left(m_1(t)|\mathscr{G}^{1}_{t}\right)}, \\
 u_{2}(t)=\frac{\mathbb{E}\left(b_2(t)x_2(t)|\mathscr{G}^{2}_{t}\right)}{\mathbb E\left(m_2(t)|\mathscr{G}^{2}_{t}\right)}+\frac{\mathbb{E}\left(m_2(t)n_{2}(t)|\mathscr{G}^{2}_{t}\right)}{\mathbb E\left(m_2(t)|\mathscr{G}^{2}_{t}\right)},
\end{aligned}
\right.
\end{equation}
where $\big((y, z_{1}, z_{2}), x_{1}, x_{2}\big)$ is a solution of the FBSDE
\begin{subequations}\label{Eq58}
\begin{numcases}{}
\nonumber
-dy(t)=\left[a(t)y(t)+b_{1}(t)\frac{\mathbb{E}\left(b_1(t)x_{1}(t)|\mathscr{G}^{1}_{t}\right)}{\mathbb E\left(m_1(t)|\mathscr{G}^{1}_{t}\right)}+b_{2}(t)\frac{\mathbb{E}\left(b_2(t)x_2(t)|\mathscr{G}^{2}_{t}\right)}{\mathbb E\left(m_2(t)|\mathscr{G}^{2}_{t}\right)}+\sum_{j=1}^{2}f_{j}(t)z_{j}(t)\right.\\
          \hspace{1.3cm} \left. +b_{1}(t)\frac{\mathbb{E}\left(m_1(t)n_{1}(t)|\mathscr{G}^{1}_{t}\right)}{\mathbb E\left(m_1(t)|\mathscr{G}^{1}_{t}\right)}+b_{2}(t)\frac{\mathbb{E}\left(m_2(t)n_{2}(t)|\mathscr{G}^{2}_{t}\right)}{\mathbb E\left(m_2(t)|\mathscr{G}^{2}_{t}\right)}+c(t)\right]dt\nonumber\\
 \label{Eq58a}                  \hspace{1.3cm}       -z_{1}(t)dw_{1}(t)-z_{2}(t)dw_{2}(t), \\
\label{Eq58b}    dx_{1}(t)= \big[a(t)x_{1}(t)-l_{1}(t)(y(t)-k_{1}(t))\big]dt+f_{1}(t)x_{1}(t)dw_{1}(t)+f_{2}(t)x_{1}(t)dw_{2}(t), \\
\label{Eq58c}    dx_{2}(t)= \big[a(t)x_{2}(t)-l_{2}(t)(y(t)-k_{1}(t))\big]dt+f_{1}(t)x_{2}(t)dw_{1}(t)+f_{2}(t)x_{2}(t)dw_{2}(t), \\
                 y(T) = \xi,\quad x_{1}(0)=-r_{1}(y(0)-h_{1}),\quad x_{2}(0)=-r_{2}(y(0)-h_{2}).
\end{numcases}
\end{subequations}
\end{proposition}


Note that since (\ref{Eq58a}) contains the conditional expectation of $x_i(t)$ with respect to $\mathscr G^i_t$ $(i=1, 2)$, (\ref{Eq58}) is new in both FBSDE and filter theories. Due to the complexity of $\mathscr G^i_t$ $(i=1, 2)$,
we are uncertain whether (\ref{Eq58}) admits a unique solution except for some special cases.


\section{Three special cases}\label{SC}   

This section focuses on solving Problem (AI) with Markovian setting, i.e., all coefficients in \eqref{Eq1} and \eqref{Eq2}
are deterministic. For the information structures (i)-(iii), we obtain the feedback Nash equilibrium points by the Riccati equations and filters of BSDEs.

\subsection{Special symmetric information: $\mathscr{G}^{1}_{t}=\mathscr{G}^{2}_{t}=\mathscr{F}^{w_{2}}_{t}$}\label{SC1}

With this symmetric information structure, we derive an explicit form of the feedback Nash equilibrium point of Problem (AI), which provides an important result for solving the asymmetric information cases in Section 3.2. Note that this result is not discussed in literature, say, \cite{WY2012}. That is why we study this case again.

Introduce two ordinary differential equations (ODEs, for short)
\begin{subequations}\label{Eq18a-b}
\begin{numcases}{}
  \dot{\alpha}_{1}-b_{1}^{2}m_{1}^{-1}\alpha_{1}^{2}-(2a+f_{2}^{2})\alpha_{1}-b_{2}^{2}m_{2}^{-1}\alpha_{1}\alpha_{2}+l_{1}=0, \label{Eq10}\\
  \dot{\beta}_{1}-(a+b_{1}^{2}m_{1}^{-1}\alpha_{1}+f_{2}^{2})\beta_{1}-b_{2}^{2}m_{2}^{-1}\alpha_{1}\beta_{2}-(b_{1}n_{1}+b_{2}n_{2}+c)\alpha_{1}-l_{1}k_{1}=0, \label{Eq11}\\
  \alpha_{1}(0)=-r_{1},\;\beta_{1}(0)=r_{1}h_{1}
\end{numcases}
\end{subequations}
and
\begin{subequations}\label{Eq21a-b}
\begin{numcases}{}
  \dot{\alpha}_{2}-b_{2}^{2}m_{2}^{-1}\alpha_{2}^{2}-(2a+f_{2}^{2})\alpha_{2}-b_{1}^{2}m_{1}^{-1}\alpha_{1}\alpha_{2}+l_{2}=0, \label{Eq15}\\
  \dot{\beta}_{2}-(a+b_{2}^{2}m_{2}^{-1}\alpha_{2}+f_{2}^{2})\beta_{2}-b_{1}^{2}m_{1}^{-1}\alpha_{2}\beta_{1}-(b_{1}n_{1}+b_{2}n_{2}+c)\alpha_{2}-l_{2}k_{2}=0, \label{Eq16}\\
  \alpha_{2}(0)=-r_{2},\;\beta_{2}(0)=r_{2}h_{2},
\end{numcases}
\end{subequations}
which will be derived step by step in Theorem \ref{Theorem3.12}. Here we omit the time variable $t$ in \eqref{Eq10},  \eqref{Eq11}, \eqref{Eq15}
and \eqref{Eq16} for simplicity.  Similar convention will be taken for the subsequent ODEs, SDEs, BSDEs and FBSDEs except for the
initial or terminal conditions.

Throughout Section \ref{SC}, we always assume that
\begin{description}
  \item[(A1).] $b_{1}^{2}(t)m_{1}^{-1}(t)=b_{2}^{2}(t)m_{2}^{-1}(t)$ and $f_{1}(t)=0, $  $t\in[0, T].$
\end{description}

The assumption provides a sufficient condition for the existence and uniqueness of solutions to \eqref{Eq18a-b} and \eqref{Eq21a-b}.

\begin{lemma}\label{Lemma3.5}
Under (A1), there exists a unique solution $(\alpha_{1}, \beta_{1}, \alpha_{2}, \beta_{2})$ to \eqref{Eq18a-b} and  \eqref{Eq21a-b}.
\end{lemma}


\emph{Proof.}\;  Let $\alpha=\alpha_{1}+\alpha_{2}$. It follows from (A1) that
\begin{equation}\label{Eq17}
  \dot{\alpha}-b_{1}^{2}m_{1}^{-1}\alpha^{2}-(2a+f_{2}^{2})\alpha+l_{1}+l_{2}=0\quad \text{on} \; (0, T], \quad \alpha(0)=-(r_{1}+r_{2}).
\end{equation}
Since \eqref{Eq17} is a standard Riccati equation, it has a unique solution $\alpha(\cdot)$.
Introduce two auxiliary equations
\begin{align}
  &\dot{\bar{\alpha}}_{1}+\left[(2a+f_{2}^{2})-b_{1}^{2}m_{1}^{-1}\alpha\right]\bar{\alpha}_{1}+l_{1}=0\quad \hbox{on} \; (0, T], \quad     \dot{\bar{\alpha}}_{1}(0)=-r_{1}, \label{Eq18}\\
  &\dot{\bar{\alpha}}_{2}+\left[(2a+f_{2}^{2})-b_{2}^{2}m_{2}^{-1}\alpha\right]\bar{\alpha}_{2}+l_{2}=0 \quad \hbox{on} \; (0, T], \quad    \dot{\bar{\alpha}}_{2}(0)=-r_{2},  \label{Eq19}
\end{align}
where $\alpha$ is the solution to \eqref{Eq17}. Obviously, \eqref{Eq18} and \eqref{Eq19} have unique solutions $\bar{\alpha}_{1}$ and $\bar{\alpha}_{2}$,
respectively. In addition, we can check that $\alpha_{1}$ and $\alpha_{2}$ in \eqref{Eq10} and \eqref{Eq15} are also the solutions to \eqref{Eq18} and \eqref{Eq19},
respectively. From the uniqueness of solution of \eqref{Eq18} with \eqref{Eq19}, it follows that
\begin{equation*}
  \bar{\alpha}_{1}=\alpha_{1},\; \bar{\alpha}_{2}=\alpha_{2},
\end{equation*}
which implies in turn that  \eqref{Eq10} and \eqref{Eq15} have the unique solutions $\alpha_{1}$ and $\alpha_{2}$, respectively.

Let $\beta=\beta_{1}+\beta_{2}$ and $\beta(0)=r_{1}h_{1}+r_{2}h_{2}$. We have
\begin{equation}\label{Eq20}
 \dot{\beta}-(a+b_{1}^{2}m_{1}^{-1}\alpha+f_{2}^{2})\beta-(b_{1}n_{1}+b_{2}n_{2}+c)\alpha-l_{1}k_{1}-l_{2}k_{2}=0\quad \hbox{on} \; (0, T],
\end{equation}
where  $\alpha$ is the solution to \eqref{Eq17}. Note that \eqref{Eq20} has a unique solution $\beta$.
Introduce
\begin{align}
   & \dot{\bar{\beta_{1}}}-(a+f_{2}^{2})\bar{\beta_{1}}-b_{2}^{2}m_{2}^{-1}\alpha_{1}\beta-(b_{1}n_{1}+b_{2}n_{2}+c)\alpha_{1}-l_{1}k_{1}=0\quad \hbox{on}\; (0, T] \label{Eq21}
 \end{align}
 with $\bar{\beta_{1}}(0)=r_{1}h_{1}$ and
 \begin{align}
   & \dot{\bar{\beta_{2}}}-(a+f_{2}^{2})\bar{\beta_{2}}-b_{1}^{2}m_{1}^{-1}\alpha_{2}\beta-(b_{1}n_{1}+b_{2}n_{2}+c)\alpha_{2}-l_{2}k_{2}=0\quad \hbox{on}\; (0, T]\label{Eq22}
\end{align}
with $\bar{\beta_{2}}(0)=r_{2}h_{2}$,
where $\alpha_{1}, \alpha_{2}$ and $ \beta$ are the solutions to \eqref{Eq18},  \eqref{Eq19} and \eqref{Eq20}, respectively. Similarly, we can prove that
\eqref{Eq11} and  \eqref{Eq16} also have unique solutions $\beta_{1}$ and $\beta_{2}$ satisfying
\begin{equation*}
  \bar{\beta_{1}}=\beta_{1},\; \bar{\beta_{2}}=\beta_{2}.
\end{equation*}
Based on the arguments above, we can derive the unique analytical expressions for $\alpha_{1}$, $\alpha_{2}$, $\beta_{1}$, $\beta_{2}$, $\alpha$ and $\beta$.
Then the proof is completed. \hfill $\Box$


\begin{theorem}\label{Theorem3.12}
 Under (A1), Problem (AI) has a unique Nash equilibrium point
 \begin{equation}\label{Eq78}
 \left\{
 \begin{aligned}
 u_{1}(t)=m_{1}^{-1}(t)b_{1}(t)\left(\alpha_{1}(t)\tilde{y}(t)+\beta_{1}(t)\right)+n_{1}(t), \\
 u_{2}(t)=m_{2}^{-1}(t)b_{2}(t)\left(\alpha_{2}(t)\tilde{y}(t)+\beta_{2}(t)\right)+n_{2}(t),
\end{aligned}
\right.
\end{equation}
where $\tilde{h}(t)=\mathbb{E}\left(h(t)\big|\mathscr{F}_{t}^{w_{2}}\right)$, $\alpha_{i}$, $\beta_{i}\, (i=1, 2)$ and $\tilde{y}$
satisfy \eqref{Eq18a-b}, \eqref{Eq21a-b} and \eqref{Eq77}, respectively.
\end{theorem}

\emph{Proof:} \noindent \textbf{(i)} We first prove that the Nash equilibrium point $(u_{1}, u_{2})$ is uniquely determined by
\begin{equation}\label{Eq41'}
 \left\{
 \begin{aligned}
 u_{1}(t)=m_{1}^{-1}(t)b_{1}(t)\tilde{x}_{1}(t)+n_{1}(t), \\
 u_{2}(t)=m_{2}^{-1}(t)b_{2}(t)\tilde{x}_{2}(t)+n_{2}(t),
\end{aligned}
\right.
\end{equation}
where $\big((y, z_{1}, z_{2}), x_{1}, x_{2}\big)$ is the solution of the FBSDE
\begin{subequations}\label{Eq3a-d'}
\begin{numcases}{}
\label{Eq3a'} -dy= \left(ay+b_{1}^{2}m_{1}^{-1}\tilde{x}_{1}+b_{2}^{2}m_{2}^{-1}\tilde{x}_{2}+f_{2}z_{2}+b_{1}n_{1}+b_{2}n_{2}+c\right)dt-\sum_{j=1}^{2}z_{j}dw_{j}, \\
\label{Eq3b'}    dx_{1}= \left[ax_{1}-l_{1}(y-k_{1})\right]dt+f_{2}x_{1}dw_{2}, \\
\label{Eq3c'}    dx_{2}= \left[ax_{2}-l_{2}(y-k_{2})\right]dt+f_{2}x_{2}dw_{2}, \\
                y(T) = \xi,\quad x_{1}(0)=-r_{1}(y(0)-h_{1}),\quad x_{2}(0)=-r_{2}(y(0)-h_{2}).
\end{numcases}
\end{subequations}
According to (A1) and Proposition \ref{Proposition3.3}, it is enough to prove the existence and uniqueness of \eqref{Eq3a-d'}.
The detail of the proof is divided into three steps.

\textbf{\emph{Step 1:}} Filtering equations.

Note that \eqref{Eq3a'} depends on the filter $\tilde x_i$. Then we need to compute the filter $\big(\tilde{y}, \tilde{z}_{2}, \tilde{x}_{1}, \tilde{x}_{2}\big)$ of $\big(y, z_{2}, x_{1}, x_{2}\big)$ with respect to $\mathscr{F}^{w_{2}}_{t}$. Applying Lemma 5.4 in \cite{Xiong2008} to \eqref{Eq3a-d'}, we get
\begin{subequations}\label{Eq5a-d}
\begin{numcases}{}
\label{Eq5}   -d\tilde{y}=\left(a\tilde{y}+b_{1}^{2}m_{1}^{-1}\tilde{x}_{1}+b_{2}^{2}m_{2}^{-1}\tilde{x}_{2}+f_{2}\tilde{z}_{2}+b_{1}n_{1}+b_{2}n_{2}+c\right)dt-\tilde{z}_{2}dw_{2}, \\
\label{Eq6}   d\tilde{x}_{1}=\left[a\tilde{x}_{1}-l_{1}(\tilde{y}-k_{1})\right]dt+f_{2}\tilde{x}_{1}dw_{2}, \\
\label{Eq7}   d\tilde{x}_{2}=\left[a\tilde{x}_{2}-l_{2}(\tilde{y}-k_{2})\right]dt+f_{2}\tilde{x}_{2}dw_{2}, \\
    \tilde{y}(T)= \mathbb{E}\left(\xi|\mathscr{F}^{w_{2}}_{T}\right),\, \tilde{ x}_{1}(0)=-r_{1}(\tilde{y}(0)-h_{1}),\, \tilde{ x}_{2}(0)=-r_{2}(\tilde{y}(0)-h_{2}).
\end{numcases}
\end{subequations}
Recall (\ref{Eq58a}). If $f_1(t)\neq0$, the generator of (\ref{Eq5}) has an additional term $f_1\tilde z_1$, which leads to a difficulty of proving the existence and uniqueness of solution to (\ref{Eq5}).

\textbf{\emph{Step 2:} }Existence and uniqueness of \eqref{Eq5a-d}.

Introduce an FBSDE
\begin{equation}\label{Eq24}
 \left\{
 \begin{aligned}
    -dp&=\left(ap+b_{1}^{2}m_{1}^{-1}n+f_{2}q+b_{1}n_{1}+b_{2}n_{2}+c\right)dt-qdw_{2}, \\
   dn&= \big[an-(l_{1}+l_{2})p+l_{1}k_{1}+l_{2}k_{2}\big]dt+f_{2}ndw_{2}, \\
   p(T)&=\mathbb{E}\left(\xi|\mathscr{F}^{w_{2}}_{T}\right), \quad n(0)=-(r_{1}+r_{2})p(0)+r_{1}h_{1}+r_{2}h_{2}.
\end{aligned}
\right.
\end{equation}
If $\big((\tilde{y}, \tilde{z}_{2}), \tilde{x}_{1}, \tilde{x}_{2}\big)$ is a solution to \eqref{Eq5a-d}, then $(n, p, q)$ is  a solution to \eqref{Eq24}, where we set
$$p=\tilde{y},\;q=\tilde{z}_{2}, \; n=\tilde{x}_{1}+\tilde{x}_{2}.$$
On the other hand, let $(p, q, n)$ be a solution to \eqref{Eq24}. Introduce an SDE
\begin{numcases}{}
\nonumber   dN_{1}=\left[aN_{1}-l_{1}(p-k_{1})\right]dt+f_{2}N_{1}dw_{2}, \\
\label{Eq27}   dN_{2}=\left[aN_{2}-l_{2}(p-k_{2})\right]dt+f_{2}N_{2}dw_{2}, \\
 \nonumber   N_{1}(0)=-r_{1}(p(0)-h_{1}),\, N_{2}(0)=-r_{2}(p(0)-h_{2}),
\end{numcases}
which has a unique solution $(N_{1}, N_{2})$ with $N_{1}+N_{2}=n$. Furthermore, we can check that
 $\big((p, q), N_{1}, N_{2}\big)$ is a solution to \eqref{Eq5a-d}. It implies that the existence and
uniqueness of \eqref{Eq5a-d} is equivalent to that of \eqref{Eq24}.
It is easy to check that \eqref{Eq24} has a unique solution $(p, q, n)$ (see, e.g., Theorem 2.3 in \cite{YJ2008}). So does
\eqref{Eq5a-d}.

\textbf{\emph{Step 3:}} Existence and uniqueness of \eqref{Eq3a-d'}.

Let $\big((\tilde{y}, \tilde{z}_{2}),\tilde{x}_{1}, \tilde{x}_{2}\big)$ be the unique solution to \eqref{Eq5a-d}. For the fixed $\tilde{x}_{1}$ and $\tilde{x}_{2}$,
we can prove that \eqref{Eq3a-d'} has a unique solution by some arguments similar to Step 2.

\noindent \textbf{(ii)} To get the feedback Nash equilibrium point, we have to establish the relationship between $\tilde y$ and $\tilde{x}_{i}\, (i=1, 2)$.
Noticing the terminal condition of \eqref{Eq3a-d'}, we set
\begin{equation}\label{Eq65}
  x_{i}=\alpha_{i}y+\beta_{i}
\end{equation}
with $ \alpha_{i}(0)=-r_{i}$  and $\beta_{i}(0)=r_{i}h_{i}, i=1, 2.$ Applying It\^{o}'s formula to $x_{1}$ in \eqref{Eq65} subject to \eqref{Eq3a'}, we obtain
\begin{multline}\label{Eq9}
  dx_{1} =\left[(\dot{\alpha}_{1}-a\alpha_{1})y-b_{1}^{2}m_{1}^{-1}\alpha_{1}\tilde{x}_{1}-b_{2}^{2}m_{2}^{-1}\alpha_{1}\tilde{x}_{2}-f_{2}\alpha_{1}z_{2}+\dot{\beta}_{1}-(b_{1}n_{1}+b_{2}n_{2}+c)\alpha_{1}\right]dt\\
           +\sum_{j=1}^{2}\alpha_{1}z_{j}dw_{j}.
\end{multline}
Substituting \eqref{Eq65} into \eqref{Eq3b'} and comparing the coefficients between \eqref{Eq3b'} and \eqref{Eq9}, we have
\begin{equation}\label{Eq12}
z_{1}=0,\quad z_{2}=f_{2}\alpha_{1}^{-1}x_{1}\equiv f_{2}y+f_{2}\alpha_{1}^{-1}\beta_{1},
\end{equation}
\begin{multline}\label{Eq67}
\left[\dot{\alpha}_{1}-(2a+f_{2}^{2})\alpha_{1}+l_{1}\right]y-b_{1}^{2}m_{1}^{-1}\alpha_{1}\tilde{x}_{1}
-b_{2}^{2}m_{2}^{-1}\alpha_{1}\tilde{x}_{2}+\dot{\beta}_{1}-\left(a+f_{2}^{2}\right)\beta_{1}\\
-(b_{1}n_{1}+b_{2}n_{2}+c)\alpha_{1}-l_{1}k_{1}=0.
\end{multline}
Taking $\mathbb{E}\left[\cdot|\mathscr{F}^{w_{2}}_{t}\right]$ on both sides of \eqref{Eq65}, \eqref{Eq12} and \eqref{Eq67}, it yields
\begin{equation}\label{Eq70}
  \tilde{x}_{i}=\alpha_{i}\tilde{y}+\beta_{i},\; i=1, 2,
\end{equation}
\begin{equation}\label{Eq68}
\tilde{z}_{1}=0,\quad \tilde{z}_{2}=f_{2}\alpha_{1}^{-1}\tilde{x}_{1}\equiv f_{2}\tilde{y}+f_{2}\alpha_{1}^{-1}\beta_{1}
\end{equation}
and
\begin{multline}\label{Eq69}
\left[\dot{\alpha}_{1}-(2a+f_{2}^{2})\alpha_{1}+l_{1}\right]\tilde{y}-b_{1}^{2}m_{1}^{-1}\alpha_{1}\tilde{x}_{1}
-b_{2}^{2}m_{2}^{-1}\alpha_{1}\tilde{x}_{2}+\dot{\beta}_{1}-\left(a+f_{2}^{2}\right)\beta_{1}\\
-(b_{1}n_{1}+b_{2}n_{2}+c)\alpha_{1}-l_{1}k_{1}=0.
\end{multline}
Plugging \eqref{Eq70} into \eqref{Eq69},  we derive \eqref{Eq18a-b}. Similarly, we have
\begin{equation}\label{Eq71}
z_{1}=0, \quad z_{2}=f_{2}\alpha_{2}^{-1}x_{2}\equiv f_{2}y+f_{2}\alpha_{2}^{-1}\beta_{2},
\end{equation}
\begin{multline}\label{Eq72}
\left[\dot{\alpha}_{2}-(2a+f_{2}^{2})\alpha_{2}+l_{2}\right]y-b_{1}^{2}m_{1}^{-1}\alpha_{2}\tilde{x}_{1}
-b_{2}^{2}m_{2}^{-1}\alpha_{2}\tilde{x}_{2}+\dot{\beta}_{2}-\left(a+f_{2}^{2}\right)\beta_{2}\\
-(b_{1}n_{1}+b_{2}n_{2}+c)\alpha_{2}-l_{2}k_{2}=0.
\end{multline}
Taking $\mathbb{E}\left[\cdot|\mathscr{F}^{w_{2}}_{t}\right]$ on both sides of \eqref{Eq71}  and \eqref{Eq72}, it yields
\begin{equation}\label{Eq73}
\tilde{z}_{1}=0,\quad \tilde{z}_{2}=f_{2}\alpha_{2}^{-1}\tilde{x}_{2}\equiv f_{2}\tilde{y}+f_{2}\alpha_{2}^{-1}\beta_{2}
\end{equation}
and
\begin{multline}\label{Eq74}
\left[\dot{\alpha}_{2}-(2a+f_{2}^{2})\alpha_{2}+l_{2}\right]\tilde{y}-b_{1}^{2}m_{1}^{-1}\alpha_{2}\tilde{x}_{1}
-b_{2}^{2}m_{2}^{-1}\alpha_{2}\tilde{x}_{2}+\dot{\beta}_{2}-\left(a+f_{2}^{2}\right)\beta_{2}\\
-(b_{1}n_{1}+b_{2}n_{2}+c)\alpha_{2}-l_{2}k_{2}=0,
\end{multline}
subject to \eqref{Eq70}. Plugging \eqref{Eq70} into \eqref{Eq74},  we derive \eqref{Eq21a-b}.

According to \eqref{Eq70}, \eqref{Eq5} is rewritten as
\begin{equation}\label{Eq76}
\left\{
\begin{aligned}
-d\tilde{y}=&\, \left[\left(a+b_{1}^{2}m_{1}^{-1}\alpha\right)\tilde{y}+f_{2}\tilde{z}_{2}+b_{1}^{2}m_{1}^{-1}\beta+b_{1}n_{1}+b_{2}n_{2}+c\right]dt-\tilde{z}_{2}dw_{2}, \\
\tilde{y}(T)=&\,\mathbb{E}\left(\xi|\mathscr{F}^{w_{2}}_{T}\right).
\end{aligned}
\right.
\end{equation}
Solving it, we get a unique solution
\begin{equation}\label{Eq77}
  \tilde{y}(t)= \mathbb{E}\left[\Gamma_{t}^{T}\mathbb{E}\left(\xi|\mathscr{F}^{w_{2}}_{T}\right)
          +\int_{t}^{T}\Gamma_{t}^{s}\left(b_{1}^{2}m_{1}^{-1}\beta+b_{1}n_{1}+b_{2}n_{2}+c\right)(s)ds\Big|\mathscr{F}^{w_{2}}_{t}\vphantom{\int_{t}^{T}}\right],
\end{equation}
where
$$
\Gamma_{t}^{s}=\exp\left\{ \int_{t}^{s}\left(a+b_{1}^{2}m_{1}^{-1}\alpha-\frac{1}{2}f_{2}^{2}\right)(r)dr\right. \\
\left.+\int_{t}^{s}f_{2}(r)dw_{2}(r)\right\},
$$
and
$\alpha=\alpha_{1}+\alpha_{2}$ and $\beta=\beta_{1}+\beta_{2}$ are uniquely given by \eqref{Eq18a-b} and \eqref{Eq21a-b}, respectively. \hfill$\Box$

\subsection{Special asymmetric information}\label{Section3.2}  


\subsubsection{ $\mathscr{G}^{1}_{t}=\mathscr{F}_{t}$  and $\mathscr{G}^{2}_{t}=\mathscr{F}^{w_{2}}_{t}.$}\label{SC2}      

In this case, $\mathbb{E}\left(x_{1}(t)|\mathscr{G}^{1}_{t}\right)=\mathbb{E}\left(x_{1}(t)|\mathscr{F}_{t}\right)=x_{1}(t)$ and $\mathbb{E}\left(x_{2}(t)|\mathscr{G}_{t}^{2}\right)=\mathbb{E}\left(x_{2}(t)|\mathscr{F}^{w_{2}}_{t}\right)=\tilde{x}_{2}(t)$. With the notations,
we get

\begin{theorem}\label{Theorem3.10}
Under (A1), Problem (AI) has a unique Nash equilibrium point denoted  by
\begin{equation}\label{Eq79}
\left\{
\begin{aligned}
u_{1}(t)=& m_{1}^{-1}(t)b_{1}(t)\big(\gamma_{1}(t)y(t)+\gamma_{2}(t)\tilde{y}(t)+\gamma_{3}(t)\big)+n_{1}(t), \\
u_{2}(t)=& m_{2}^{-1}(t)b_{2}(t)\big(\alpha_{2}(t)\tilde{y}(t)+\beta_{2}(t)\big)+n_{2}(t).
\end{aligned}
\right.
\end{equation}
Here $\tilde{y}$ and  $y$ are given by \eqref{Eq77} and \eqref{Eq97}, respectively; $(\alpha_{2}, \beta_{2})$ and $(\gamma_{1}, \gamma_{2}, \gamma_{3})$ are
the solutions to \eqref{Eq21a-b} and \eqref{Eq28a-c}, respectively.
\end{theorem}

Note that even through the player 1 has access to the complete information, the information available to the player 2 has an effect on the control policy of the player 1 via $\tilde y(t)$. This is an interesting phenomenon indeed.

\emph{Proof.} We complete this proof by two steps.

\emph{\textbf{Step 1:}} We prove that under (A1), the Nash equilibrium point is uniquely determined by
\begin{equation}\label{Eq60}
 \left\{
 \begin{aligned}
 u_{1}(t)=m_{1}^{-1}(t)b_{1}(t)x_{1}(t)+n_{1}(t), \\
 u_{2}(t)=m_{2}^{-1}(t)b_{2}(t)\tilde{x}_{2}(t)+n_{2}(t),
\end{aligned}
\right.
\end{equation}
where $\big((y, z_{1}, z_{2}), x_{1}, x_{2}\big)$
is the solution of the FBSDE
 \begin{subequations}\label{Eq59}
\begin{numcases}{}
\label{Eq59a}   -dy=\left(ay+b_{1}^{2}m_{1}^{-1}x_{1}+b_{2}^{2}m_{2}^{-1}\tilde{x}_{2}+f_{2}z_{2}+b_{1}n_{1}+b_{2}n_{2}+c\right)dt-\sum_{j=1}^{2}z_{j}dw_{j}, \\
\label{Eq59b}    dx_{1}= \left[ax_{1}-l_{1}(y-k_{1})\right]dt+f_{2}x_{1}dw_{2}, \\
\label{Eq59c}  dx_{2}= \left[ax_{2}-l_{2}(y-k_{2})\right]dt+f_{2}x_{2}dw_{2}, \\
\label{Eq59d}  y(T) = \xi,\quad x_{1}(0)=-r_{1}(y(0)-h_{1}), \quad  x_{2}(0)=-r_{2}(y(0)-h_{2}).
\end{numcases}
\end{subequations}
Similar to Theorem \ref{Theorem3.12}, we only need to prove the existence and uniqueness of \eqref{Eq59}.
It is easy to see that the optimal filter $\big(\tilde{y}, \tilde{z}_{2}, \tilde{x}_{1}, \tilde{x}_{2}\big)$ of $\big(y,  z_{2}, x_{1}, x_{2}\big)$ in \eqref{Eq59}
still satisfies \eqref{Eq5a-d}. Thus, $\tilde{y}$ is given by \eqref{Eq77}, and $\tilde{x}_{2}$ is uniquely represented by $\tilde{y}$ as shown in \eqref{Eq70}.
Then \eqref{Eq59a} with \eqref{Eq59b} is rewritten as
\begin{subequations}\label{Eq62}
\begin{numcases}{}
 \label{Eq62a} -dy=\left(ay+f_{2}z_{2}+b_{1}^{2}m_{1}^{-1}x_{1}+b_{2}^{2}m_{2}^{-1}\alpha_{2}\tilde{y}+b_{2}^{2}m_{2}^{-1}\beta_{2}+b_{1}n_{1}+b_{2}n_{2}+c\right)dt\nonumber\\
                   \hspace{12mm}            -\sum_{j=1}^{2}z_{j}dw_{j},\\
\label{Eq62b}  dx_{1}= \left[ax_{1}-l_{1}(y-k_{1})\right]dt+f_{2}x_{1}dw_{2}, \\
\label{Eq62c}  y(T) = \xi,\quad x_{1}(0)=-r_{1}(y(0)-h_{1}).
\end{numcases}
\end{subequations}
Thanks to Theorem 2.3 in \cite{YJ2008}, \eqref{Eq62} has a unique solution $\big(y, z_{1}, z_{2}, x_{1}\big)$. Substituting $y$ in \eqref{Eq62} into
\eqref{Eq59c} and \eqref{Eq59d}, \eqref{Eq59c} has a unique solution $x_{2}$. Therefore, \eqref{Eq59} is uniquely solvable.

\emph{\textbf{Step 2:}} We verify that the feedback Nash equilibrium
point is shown as \eqref{Eq79}. According to \eqref{Eq62a} and
\eqref{Eq62b} together with the initial condition in \eqref{Eq62c},
we set
\begin{equation}\label{Eq61}
  x_{1}=\gamma_{1}y+\gamma_{2}\tilde{y}+\gamma_{3}
\end{equation}
with $\gamma_{1}(0)=-r_{1}, \gamma_{2}(0)=0, \gamma_{3}(0)=r_{1}h_{1}.$ Applying It\^{o}'s formula to $x_{1}$ in \eqref{Eq61}, we have
\begin{multline} \label{Eq26}
dx_{1}=\Big\{\left(\dot{\gamma_{1}}-a\gamma_{1}\right)y+\left(\dot{\gamma_{2}}-(a+b_{1}^{2}m_{1}^{-1}\alpha)\gamma_{2}\right)\tilde{y}
            -b_{1}^{2}m_{1}^{-1}\gamma_{1}x_{1}-b_{2}^{2}m_{2}^{-1}\gamma_{1}\tilde{x}_{2}-\gamma_{1}f_{2}z_{2}\\
-\gamma_{2}f_{2}\tilde{z}_{2}+\dot{\gamma_{3}}-(b_{1}n_{1}+b_{2}n_{2}+c)\gamma_{1}
-(b_{1}n_{1}+b_{2}n_{2}+c+b_{1}^{2}m_{1}^{-1}\beta)\gamma_{2}\Big\}dt\\
+\gamma_{1}z_{1}dw_{1}+\left(\gamma_{1}z_{2}+\gamma_{2}\tilde{z}_{2}\right)dw_{2}
\end{multline}
with $\tilde{x}_{2}=\alpha_{2}\tilde{y}+\beta_{2}$ and $\tilde{z}_{2}=f_{2}\tilde{y}+f_{2}\alpha_{2}^{-1}\beta_{2}.$
 Comparing  \eqref{Eq62b} with \eqref{Eq26}, we get
 \begin{equation}\label{Eq32}
   z_{1}=0,\quad z_{2}=f_{2}y+f_{2}\gamma_{1}^{-1}\gamma_{3}-f_{2}\gamma_{1}^{-1}\gamma_{2}\alpha_{2}^{-1}\beta_{2},
 \end{equation}
 \begin{align}\label{Eq87}
  & \Big[\dot{\gamma}_{2}-\left(a+b_{1}^{2}m_{1}^{-1}\alpha+f_{2}^{2}+b_{1}^{2}m_{1}^{-1}\gamma_{1}\right)\gamma_{2}
    -b_{2}^{2}m_{2}^{-1}\alpha_{2}\gamma_{1}\Big]\tilde{y}+\Big[\dot{\gamma}_{1}-(a+f_{2}^{2})\gamma_{1}-b_{1}^{2}m_{1}^{-1}\gamma_{1}^{2}\Big]y\nonumber \\
  & +\dot{\gamma}_{3}-\Big(f_{2}^{2}+b_{1}^{2}m_{1}^{-1}\gamma_{1}\Big)\gamma_{3}
    -\Big(b_{1}n_{1}+b_{2}n_{2}+c
    +b_{2}^{2}m_{2}^{-1}\beta_{2}\Big)\gamma_{1}
     -\Big(b_{1}n_{1}+b_{2}n_{2}+c+b_{1}^{2}m_{1}^{-1}\beta\Big)\gamma_{2}\nonumber\\
  & =(a\gamma_{1}-l_{1})y+a\gamma_{2}\tilde{y}+a\gamma_{3}+l_{1}k_{1}.
 \end{align}
Then we have
\begin{subequations}\label{Eq28a-c}
\begin{numcases}{}
\dot{\gamma_{1}}-b_{1}^{2}m_{1}^{-1}\gamma_{1}^{2}-(2a+f_{2}^{2})\gamma_{1}+l_{1}=0,\label{Eq28}\\
\dot{\gamma}_{2}-(2a+b_{1}^{2}m_{1}^{-1}\alpha+f_{2}^{2}+b_{1}^{2}m_{1}^{-1}\gamma_{1})\gamma_{2}-b_{2}^{2}m_{2}^{-1}\alpha_{2}\gamma_{1}=0,\label{Eq29}\\
\dot{\gamma}_{3}-(a+f_{2}^{2}+b_{1}^{2}m_{1}^{-1}\gamma_{1})\gamma_{3}-l_{1}k_{1}-(b_{1}n_{1}+b_{2}n_{2}+c+b_{2}^{2}m_{2}^{-1}\beta_{2})\gamma_{1}\nonumber\\
\hspace{58mm}        -(b_{1}n_{1}+b_{2}n_{2}+c+b_{1}^{2}m_{1}^{-1}\beta)\gamma_{2}=0,\label{Eq30}\\
\gamma_{1}(0)=-r_{1},\; \gamma_{2}(0)=0,\; \gamma_{3}(0)=r_{1}h_{1},
\end{numcases}
\end{subequations}
which has a unique solution $(\gamma_{1}, \gamma_{2}, \gamma_{3}).$  Substituting \eqref{Eq61} into \eqref{Eq62a}, we derive
\begin{equation}\label{Eq97}
y(t)=\mathbb{E}\left(\xi\Upsilon_{t}^{T}+\int_{t}^{T}\Upsilon_{t}^{s}g_{2}(s)ds|\mathscr{F}_{t}\right)
\end{equation}
with
\begin{align*}
&  \Upsilon_{t}^{s}=\exp\left\{\int_{t}^{s}\left(g_{1}(r)-\frac{1}{2}f_{2}^{2}(r)\right)dr
           +\int_{t}^{s}f_{2}(r)dw_{2}(r)\right\},\\
& g_{1}=a+b_{1}^{2}m_{1}^{-1}\gamma_{1},\\
& g_{2}=(b_{1}^{2}m_{1}^{-1}\gamma_{2}+b_{2}^{2}m_{2}^{-1}\alpha_{2})\tilde{y}+b_{1}^{2}m_{1}^{-1}\gamma_{3}+b_{2}^{2}m_{2}^{-1}\beta_{2}+b_{1}n_{1}+b_{2}n_{2}+c.
\end{align*}
Then the proof is completed. \hfill$\Box$

\begin{remark}
The above arguments can also be used to solve the case of $\mathscr{G}^{1}_{t}=\mathscr{F}_{t}$ and $\mathscr{G}^{2}_{t}=\mathscr{F}^{w_1}_{t}$.
We omit it here.
\end{remark}

\subsubsection{ $\mathscr{G}^{1}_{t}=\mathscr{F}^{w_{1}}_{t}$  and $\mathscr{G}^{2}_{t}=\mathscr{F}^{w_{2}}_{t}.$}\label{SC3}

We assume that

\vspace{2mm}

\noindent \textbf{(A2).}  $f_{2}(t)=0,\, t\in [0, T]$.

\vspace{2mm}

With the assumption, the filter of $\big(y, z_{1}, z_{2}, x_{1}, x_{2}\big)$ in \eqref{Eq58} with respect to $\mathscr{F}^{w_{1}}_{t}$ is existent and unique. Then we derive the following feedback  Nash equilibrium point.

\begin{theorem}\label{Theorem3.6}
 Under (A1) and (A2), the feedback Nash equilibrium point of Problem (AI) is uniquely denoted  by
 \begin{equation}\label{Eq80}
 \left\{
 \begin{aligned}
 u_{1}(t)=m_{1}^{-1}(t)b_{1}(t)\big(\gamma_{1}(t)\hat{y}(t)+\gamma_{2}(t)\mathbb{E}y(t)+\gamma_{3}(t)\big)+n_{1}(t), \\
 u_{2}(t)=m_{2}^{-1}(t)b_{2}(t)\big(\tau_{1}(t)\tilde{y}(t)+\tau_{2}(t)\mathbb{E}y(t)+\tau_{3}(t)\big)+n_{2}(t).
\end{aligned}
\right.
\end{equation}
Here $\mathbb{E}y, \hat{y} $ and $\tilde{y}$ are given below in \eqref{Eq49}, \eqref{Eq51} and \eqref{Eq53}, respectively; $\gamma_{i}$ and $\tau_{i}\, (i=1, 2, 3)$
are uniquely determined by \eqref{Eq28a-c} and \eqref{Eq31a-c} with $f_{2}$ replaced by $0,$ respectively.
\end{theorem}

\emph{Proof:} \textbf{Firstly}, we prove under (A1) and (A2), Problem (AI) has a unique Nash equilibrium point determined  by
\begin{equation}\label{Eq42}
 \left\{
 \begin{aligned}
 u_{1}(t)=m_{1}^{-1}(t)b_{1}(t)\hat{x}_{1}(t)+n_{1}(t), \\
 u_{2}(t)=m_{2}^{-1}(t)b_{2}(t)\tilde{x}_{2}(t)+n_{2}(t),
\end{aligned}
\right.
\end{equation}
where  $\big((y, z_{1}, z_{2}), x_{1}, x_{2}\big)$ is the solution of the FBSDE
 \begin{subequations}\label{Eq43}
\begin{numcases}{}
\label{Eq43a}   -dy=\Big[ay+b_{1}^{2}m_{1}^{-1}\hat{x}_{1}+b_{2}^{2}m_{2}^{-1}\tilde{x}_{2}+b_{1}n_{1}+b_{2}n_{2}+c\Big]dt-\sum_{j=1}^{2}z_{j}dw_{j}, \\
\label{Eq43b}    dx_{1}= \big[ax_{1}-l_{1}(y-k_{1})\big]dt, \\
\label{Eq43c}  dx_{2}= \big[ax_{2}-l_{2}(y-k_{2})\big]dt, \\
  y(T) = \xi,\quad x_{1}(0)=-r_{1}(y(0)-h_{1}),\quad x_{2}(0)=-r_{2}(y(0)-h_{2}).
\end{numcases}
\end{subequations}

\vspace{2mm}

Once again, it is enough to prove the existence and uniqueness of the solution to \eqref{Eq43}. By the method similar to that of Theorem \ref{Theorem3.12},
the optimal filters $\hat{y}$ and $\hat{x}_{1}$ of $y$ and $x_{1}$ in \eqref{Eq43a} and \eqref{Eq43b} with respect to $\mathscr{F}^{w_{1}}_{t}$ are governed by
\begin{subequations}\label{Eq44}
\begin{numcases}{}
\label{Eq44a} -d\hat{y}=\Big[a\hat{y}+b_{1}^{2}m_{1}^{-1}\hat{x}_{1}+b_{2}^{2}m_{2}^{-1}\mathbb{E}x_{2}+b_{1}n_{1}+b_{2}n_{2}+c\Big]dt
                          -\hat{z}_{1}dw_{1}, \\
\label{Eq44b}   d\hat{x}_{1}= \big[a\hat{x}_{1}-l_{1}(\hat{y}-k_{1})\big]dt, \\
\label{Eq44c}   \hat{y}(T) = \mathbb{E}\left(\xi|\mathscr{F}^{w_{1}}_{T}\right),\quad \hat{x}_{1}(0)=-r_{1}(\hat{y}(0)-h_{1}).
\end{numcases}
\end{subequations}
Here $\mathbb{E}\eta$ stands for the expectation $\mathbb{E}\big(\eta(t)\big)$ of $\eta(t)$.
Similarly, we obtain the optimal filters $\tilde{y}$ and $\tilde{x}_{2}$ of $y$ and $x_{2}$, in \eqref{Eq43a} and \eqref{Eq43c}, with respect to
$\mathcal{F}^{w_{2}}_{t}$ as follows:
\begin{subequations}\label{Eq45}
\begin{numcases}{}
\label{Eq45a}-d\tilde{y}=\Big[a\tilde{y}+b_{1}^{2}m_{1}^{-1}\mathbb{E}x_{1}+b_{2}^{2}m_{2}^{-1}\tilde{x}_{2}+b_{1}n_{1}+b_{2}n_{2}+c\Big]dt
                   -\tilde{z}_{2}dw_{2}, \\
\label{Eq45b}  d\tilde{x}_{2}= \big[a\tilde{x}_{2}-l_{2}(\tilde{y}-k_{2})\big]dt, \\
\label{Eq45c}   \tilde{y}(T) =  \mathbb{E}\left(\xi|\mathscr{F}^{w_{2}}_{T}\right),\quad  \tilde{ x}_{2}(0)=-r_{2}(\tilde{y}(0)-h_{2}).
\end{numcases}
\end{subequations}
On the other hand,  $\mathbb{E}x_{1}$ and  $\mathbb{E}x_{2}$  together with  $\mathbb{E}y$  satisfy an ordinary differential equation
\begin{subequations}\label{Eq46}
\begin{numcases}{}
\label{Eq46a}    -\dot{\mathbb{E}}y=a\mathbb{E}y+b_{1}^{2}m_{1}^{-1}\mathbb{E}x_{1}+b_{2}^{2}m_{2}^{-1}\mathbb{E}x_{2}+b_{1}n_{1}+b_{2}n_{2}+c, \\
\label{Eq46b}    \dot{\mathbb{E}}x_{1}=a\mathbb{E}x_{1}-l_{1}\mathbb{E}y+l_{1}k_{1}, \\
\label{Eq46c}      \dot{\mathbb{E}}x_{2}=a\mathbb{E}x_{2}-l_{2}\mathbb{E}y+l_{2}k_{2}, \\
     \mathbb{E}y(T) = \mathbb{E}\xi,\quad \mathbb{E}x_{1}(0)=-r_{1}(\mathbb{E}y(0)-h_{1}),\quad \mathbb{E}x_{2}(0)=-r_{2}(\mathbb{E}y(0)-h_{2}),
\end{numcases}
\end{subequations}
where $\dot{\mathbb{E}}\eta$ denotes $\frac{d\mathbb{E}(\eta(t))}{dt}$ for $\eta=y,  x_{1},  x_{2}.$ Using the method shown in Step 2 of
Theorem \ref{Theorem3.12} again, we conclude that \eqref{Eq46} has a unique solution $(\mathbb{E}y, \mathbb{E}x_{1},$ $\mathbb{E}x_{2})$
under (A1) and (A2) (see the diffusion degenerate case of Theorem 2.3 in Yu and Ji \cite{YJ2008}). Plugging $\mathbb{E}x_{2}$ and $\mathbb{E}x_{1}$ into \eqref{Eq44}
and \eqref{Eq45}, we conclude that \eqref{Eq44} and \eqref{Eq45} have the unique solutions $\big((\hat{y}, \hat{z}_{1}), \hat{x}_{1}\big)$ and $\big((\tilde{y}, \tilde{z}_{2}), \tilde{x}_{2}\big)$, respectively. For the fixed $\hat{x}_{1}$ and $\tilde{x}_{2}$,  \eqref{Eq43} is decoupled, then it has a unique solution $(y, z_{1}, z_{2}, x_{1}, x_{2})$.

\vspace{2mm}

\textbf{Subsequently}, we verify that \eqref{Eq80} is the feedback
Nash equilibrium point. Since the required calculuses are similar to
those of Sections \ref{SC1} and \ref{SC2}, we omit unnecessary
technical details, but present key steps for the convenience of the
reader.

\vspace{2mm}

The relationship between $\mathbb{E}x_{i}$ and $\mathbb{E}y$ is
\begin{equation}\label{Eq48}
  \mathbb{E}x_{i}=\alpha_{i}\mathbb{E}y+\beta_{i} \qquad (i=1,2),
\end{equation}
where
$\alpha_{i}$, $\beta_{i}$, $\alpha$ and $\beta$ are the unique solutions to \eqref{Eq18a-b}-\eqref{Eq17} and \eqref{Eq20} with $f_{i}(\cdot)=0\; (i=1,2)$, and
\begin{equation}\label{Eq49}
  \mathbb{E}y(t)=\bar{\Gamma}_{t}^{T}\mathbb{E}\xi+\int_{t}^{T}\bar{\Gamma}_{t}^{s}\Big[\big(b_{1}^{2}(s)m_{1}^{-1}(s)\beta(s)+b_{1}n_{1}+b_{2}n_{2}+c(s)\big)\Big]ds
\end{equation}
with $$\bar{\Gamma}_{t}^{s}=\exp\left\{\int_{t}^{s}\big[a(r)+b_{1}^{2}(r)m_{1}^{-1}(r)\alpha(r)\big]dr\right\}.$$

The filter $\hat{x}_{1}$ is written as
\begin{equation}\label{Eq50}
  \hat{x}_{1}=\gamma_{1}\hat{y}+\gamma_{2}\mathbb{E}y+\gamma_{3},
\end{equation}
where
$\gamma_{i}\, (i=1, 2, 3)$ is the solution to \eqref{Eq28a-c}  with $f_{i}(\cdot)=0\; (i=1,2)$, and
\begin{equation}\label{Eq51}
  \hat{y}(t)=\Xi_{t}^{T} \mathbb{E}\left(\xi|\mathscr{F}^{w_{1}}_{t}\right)+\int_{t}^{T}\Xi_{t}^{s}g_{3}(s)ds
\end{equation}
with $$\Xi_{s}^{t}=\exp\left\{\int_{t}^{s}\big[a(r)+b_{1}^{2}(r)m_{1}^{-1}(r)\gamma_{1}(r)\big]dr \right\}$$
and $$g_{3}=\big(b_{2}^{2}m_{2}^{-1}\alpha_{2}+b_{1}^{2}m_{1}^{-1}\gamma_{2}\big)\mathbb{E}y+b_{1}^{2}m_{1}^{-1}\gamma_{3}+b_{2}^{2}m_{2}^{-1}\beta_{2}+b_{1}n_{1}+b_{2}n_{2}+c.$$
Also, $\tilde{x}_{2}$ is written as
\begin{equation}\label{Eq52}
  \tilde{x}_{2}=\tau_{1}\tilde{y}+\tau_{2}\mathbb{E}y+\tau_{3},
\end{equation}
where $(\tau_{1}, \tau_{2}, \tau_{3})$ is the unique solution to
\begin{subequations}\label{Eq31a-c}
\begin{numcases}{}
   \dot{\tau_{1}}-b_{2}^{2}m_{2}^{-1}\tau_{1}^{2}-2a\tau_{1}+l_{2}=0,\label{Eq31a}\\
    \dot{\tau_{2}}-\big(2a+b_{1}^{2}m_{1}^{-1}\alpha+b_{2}^{2}m_{2}^{-1}\tau_{1}\big)\tau_{2}
                  -b_{1}^{2}m_{1}^{-1}\alpha_{1}\tau_{1}=0,\label{Eq31b}\\
   \dot{\tau_{3}}-(a+b_{2}^{2}m_{2}^{-1}\tau_{1})\tau_{3}-(b_{1}n_{1}+b_{2}n_{2}+c+ b_{1}^{2}m_{1}^{-1}\beta_{1})\tau_{1}\nonumber\\
   \hspace{37mm}-(b_{1}n_{1}+b_{2}n_{2}+c+ b_{1}^{2}m_{1}^{-1}\beta)\tau_{2}-l_{2}k_{2}=0,\label{Eq31c}\\
   \tau_{1}(0)=-r_{2}, \; \tau_{2}(0)=0, \; \tau_{3}(0)=r_{2}h_{2}.
\end{numcases}
\end{subequations}
Then we derive
\begin{equation}\label{Eq53}
  \tilde{y}(t)=\Psi_{t}^{T} \mathbb{E}\left(\xi|\mathscr{F}^{w_{2}}_{t}\right)+\int_{t}^{T}\Psi_{t}^{s}g_{4}(s)ds
\end{equation}
with $$\Psi_{t}^{s}=\exp\left\{\int_{t}^{s}\big[a(r)+b_{2}^{2}(r)m_{2}^{-1}(r)\tau_{1}(r)\big]dr \right\}$$ and
$$g_{4}=\big(b_{2}^{2}m_{2}^{-1}\tau_{2}+b_{1}^{2}m_{1}^{-1}\alpha_{1}\big)\mathbb{E}y+b_{1}^{2}m_{1}^{-1}\beta_{1}+b_{2}^{2}m_{2}^{-1}\tau_{3}
          +b_{1}n_{1}+b_{2}n_{2}+c.$$
Thus, \eqref{Eq80} is the feedback Nash equilibrium point. Then the
proof is completed. \hfill$\Box$


\subsection{Numerical example}

 This section is devoted to illustrating the above results by a numerical example. Without loss of generality, we let $a=1, b_{1}=1, b_{2}=2, f_{1}=0, f_{2}=1, l_{1}=2,  l_{2}=4,  m_{1}=1,  m_{2}=4,  n_{1}=n_{2}=k_{1}=k_{2}=c=0, r_{1}=2, r_{2}=1, h_{1}=h_{2}=0$ in Section 3.2.1. Solving \eqref{Eq21a-b}, \eqref{Eq77}, (\ref{Eq28a-c}) and (\ref{Eq97}), we get
\begin{equation*}
 \left\{
 \begin{aligned}
&  \alpha_{2}(t)=\frac{24e^{5t}-20e^{4t}+1}{1-6e^{5t}}, \\
&  \beta_{2}(t)=0,\\
& \tilde{y}(t)=\exp\left\{\frac{3}{2}(T-t)+\ln\frac{6e^{5t}-1}{6e^{5T}-1}\right\}\\
& \hspace{11mm} \times\mathbb{E}\left(\xi e^{w_{2}(T)-w_{2}(t)}|\mathscr{F}^{w_{2}}_{t}\right),\\
&  \gamma_{1}(t)=\frac{2e^{-3t}+4}{e^{-3t}-4}, \\
&  \gamma_{2}(t)=1+\frac{54e^{5t}-40te^{4t}-\frac{80}{3}e^{7t}-\frac{245}{6}e^{4t}-\frac{3}{2}}{24e^{8t}-6e^{5t}-4e^{3t}+1},\\
&  \gamma_{3}(t)=0, \\
&  y(t)=\mathbb{E}\left[\xi\Upsilon_{t}^{T}+\int_{t}^{T}\Upsilon_{t}^{s}\left(\gamma_{2}(s)+\alpha_{2}(s)\right) \tilde{y}(s)ds|\mathscr{F}_{t}\right]
\end{aligned}
\right.
\end{equation*}
with $\Upsilon_{t}^{s}= \exp\left\{\frac{1}{2}(s-t)+\ln\frac{4e^{3t}-1}{4e^{3s}-1}+w_{2}(s)-w_{2}(t)\right\}$. Then Theorem 3.2 implies that the feedback Nash equilibrium point is uniquely denoted by
 \begin{equation*}
 \left\{
 \begin{aligned}
 &u_{1}(t)=\gamma_{1}(t)y(t)+\gamma_{2}(t)\tilde{y}(t), \\
 &u_{2}(t)=\frac{1}{2}\alpha_{2}(t)\tilde{y}(t).
\end{aligned}
\right.
\end{equation*}
Similarly, we can also perform numerical computations of the Nash equilibrium points in Theorems 3.1 and 3.3. We omit them for simplicity.


\section{Concluding remarks}\label{CR}  

This paper studies an LQ non-zero sum differential game problem, where the information available to the players is asymmetric, and the game system is a BSDE
rather than an SDE. Using the filters of FBSDEs and the existence and uniqueness of FBSDEs, we obtain the feedback Nash equilibrium points of the game problem
with observable information generated by Brownian motions. Also, we prove the uniqueness of the equilibrium points.

Three observable filtrations (see the information structures (i)-(iii) in Section \ref{PF}) are described to  classify the information available to the two players. Although the observable information of the player 2 is same in these three cases, the control policy of the player 2 varies according to the control policy of the player 1. This interesting phenomenon reflects the game behavior of these two players very nicely.
The results in Section \ref{SC} are based on $f_{1}(t)=f_{2}(t)=0$. If $f_{1}(t)f_{2}(t)\neq0$, it is difficult to prove the existence and uniqueness of the Nash equilibrium
point. We shall come back to this case in a future work.






\section*{Appendix}

\vspace{2mm} In this appendix, we use a few novel examples to illustrate the reasonability and significance of studying the special cases in Section \ref{SC}.

\vspace{2mm} \textbf{Example A.1.} Consider a controlled BSDE
\begin{equation}\label{mEq8}
\left\{
  \begin{aligned}
    -dy^{v_{1}, v_{2}}&(t) = g\big(t, y^{v_{1}, v_{2}}(t), z_{1}^{v_{1}, v_{2}}(t), z_{2}^{v_{1}, v_{2}}(t),\\
    &\hspace{12mm}  v_{1}(t), v_{2}(t)\big)dt\\
            &\ -z_{1}^{v_{1}, v_{2}}(t)dw_{1}(t)-z_{2}^{v_{1}, v_{2}}(t)dw_{2}(t), \\
   y^{v_{1}, v_{2}}& (T)= x(T)
  \end{aligned}
\right.
\end{equation}
with
$$
\left\{\begin{aligned}
dx(t)= &\ b\big(t, x(t)\big)dt+\delta_{1}(t)dw_{1}(t)+\delta_{2}(t)dw_{2}(t),\\
x(0)= &\ 0.
\end{aligned}
\right.
$$
Here $\delta_{1}$ and $\delta_{2}$ are uniformly bounded and
deterministic; $b$ and $g$ are deterministic and satisfy certain
conditions which guarantee the existence and uniqueness of solution
to \eqref{mEq8}; and $v_{1}$ and $v_{2}$ are control processes for
the player 1 and the player 2, respectively. Note that $x$ is not
controlled, and $y^{v_{1}, v_{2}}$ is coupled with $x$ at the
terminal time $T$. Cost functional for the player $i$ $(i=1, 2)$ is
of the form
\begin{multline}\label{mEq2}
 J_{i}\big(v_{1}(\cdot), v_{2}(\cdot)\big)=\frac{1}{2} \mathbb{E}\left\{ \int_{0}^{T} l_{i}\big(t, y^{v_{1}, v_{2}}(t), z_{1}^{v_{1}, v_{2}}(t),z_{2}^{v_{1}, v_{2}}(t), v_{1}(t), v_{2}(t)\big)dt+r_{i}\big(y^{v_{1}, v_{2}}(0)\big)\vphantom{\int_{0}^{T}}\right\},
\end{multline}
where $l_{i}$ and $r_{i}$ are deterministic, and satisfy certain
integrability conditions. Assume that the player 1 has access to the
complete information $\mathscr{F}_{t}$, i.e., the player 1 selects
his/her control process $v_{1}$ according to $\mathscr{F}_{t}$.
However, the player 2 can only partially observe the state $(x,
y^{v_{1}, v_{2}}, z_{1}^{v_{1}, v_{2}}, z_{2}^{v_{1}, v_{2}})$
through a noisy process
\begin{equation}\label{mEq6}
\left\{\begin{aligned}
  dW_{2}(t)=&\ h\big(t, x(t)\big)dt+dw_{2}(t), \\
  W_{2}(0)=&\ 0,
  \end{aligned}
  \right.
\end{equation}
where $h$ is deterministic and uniformly bounded. Define the
admissible control sets
$$\mathscr{X}_{1}=\left\{v_{1}(\cdot); v_{1}(t)\  \mbox{is}\  \mathscr{F}_{t}\mbox{-adapted and square-integrable}\right\}$$
and
\begin{equation*}
  \mathscr{X}_{2}=\left\{v_{2}(\cdot); v_{2}(t)\  \mbox{is}\   \mathscr{F}_{t}^{W_{2}}\mbox{-}\mbox{adapted and } \mbox{ square-integrable}\vphantom{\mathscr{F}_{t}^{W_{2}}}\right\}.
\end{equation*}
Then the game problem is stated as follows.

\textbf{Problem (A.1)}. Find a pair of admissible controls $(u_{1}, u_{2})$ such that
\begin{equation*}
 \left\{
 \begin{aligned}
& J_{1}(u_{1}(\cdot), u_{2}(\cdot)) = \min\limits_{v_{1}(\cdot)\in \mathscr{X}_{1}} J_{1}(v_{1}(\cdot), u_{2}(\cdot)),\\
& J_{2}(u_{1}(\cdot), u_{2}(\cdot)) = \min\limits_{v_{2}(\cdot)\in \mathscr{X}_{2}} J_{2}(u_{1}(\cdot), v_{2}(\cdot)),
\end{aligned}
 \right.
\end{equation*}
subject to \eqref{mEq8}, \eqref{mEq2} and \eqref{mEq6}.

In the sequel, we wish to simplify Problem (A.1) by an equivalent transformation. Set
\begin{equation*}
  \rho_{1}(t)=\exp\left\{-\int_{0}^{t}h(s, x(s)) dw_{2}(s)-\frac{1}{2}\int_{0}^{t}|h(s, x(s))|^{2}ds\right\}
\end{equation*}
and $$\frac{d\mathbb{Q}}{d\mathbb{P}}\Big |_{\mathscr{F}_{T}}=\rho_{1}(T).$$
Since $h$ is bounded, Girsanov theorem implies that $\mathbb{Q}$ is
a new probability measure, and thus $(w_{1}, W_{2})$ is a standard
Brownian motion under $\mathbb{Q}$. Plugging (\ref{mEq6}) into
(\ref{mEq8}), we have
\begin{equation}\label{mEq10}
\left\{
  \begin{aligned}
  -dy^{v_{1}, v_{2}}(t) = &\ \big[g\big(t, y^{v_{1}, v_{2}}(t), z_{1}^{v_{1}, v_{2}}(t), z_{2}^{v_{1}, v_{2}}(t), v_{1}(t), v_{2}(t)\big)+h(t, x(t))z_{2}^{v_{1}, v_{2}}(t)\big]dt\\
            &\ -z_{1}^{v_{1}, v_{2}}(t)dw_{1}(t)-z_{2}^{v_{1}, v_{2}}(t)dW_{2}(t), \\
y^{v_{1}, v_{2}}(T) = &\ x(T)
  \end{aligned}
\right.
\end{equation}
with
$$
\left\{
  \begin{aligned}
     dx(t) = &\ \big[b\big(t, x(t)\big)-\delta_{2}(t)h(t, x(t))\big]dt+\delta_{1}(t)dw_{1}(t)+\delta_{2}(t)dW_{2}(t),\\
     x(0)= &\ 0.
  \end{aligned}
\right.
$$
On the other hand,
\begin{equation*}
  \rho_{1}^{-1}(t)=\exp\left\{\int_{0}^{t}h(s, x(s)) dW_{2}(s)-\frac{1}{2}\int_{0}^{t}|h(s, x(s))|^{2}ds\right\}.
\end{equation*}
Then \eqref{mEq2} is rewritten as
\begin{multline}\label{mEq11}
  \breve{J}_{i}\big(v_{1}(\cdot), v_{2}(\cdot)\big)=\frac{1}{2} \mathbb{E}_{\mathbb{Q}}\left\{ \int_{0}^{T}\rho_{1}^{-1}(t)l_{i}\big(t, y^{v_{1}, v_{2}}(t), z_{1}^{v_{1}, v_{2}}(t), z_{2}^{v_{1}, v_{2}}(t), v_{1}(t), v_{2}(t)\big)dt\right.\\
\left.+r_{i}\big(y^{v_{1}, v_{2}}(0)\big)\vphantom{\int_{0}^{T}}\right\}.
\end{multline}
We can check that $\mathscr{F}_{t}=\mathscr{F}^{w_{1}, W_{2}}_{t}$.
So $\mathscr{X}_{1}$ is equivalent to the  admissible control set
$$
  \mathscr{Y}_{1}=\left\{v_{1}(\cdot); v_{1}(t)\; \hbox{is an}\; \mathscr{F}^{w_{1}, W_{2}}_{t}\hbox{-adapted and}  \hbox{ square-integrable process}\right\}.
$$
Now Problem (A.1) can be equivalently stated as follows.

\textbf{Problem (A.1')}. Find a pair of admissible controls $(u_{1}, u_{2})$ such that
\begin{equation*}
 \left\{
 \begin{aligned}
& \breve{J}_{1}(u_{1}(\cdot), u_{2}(\cdot)) = \min\limits_{v_{1}(\cdot)\in \mathscr{Y}_{1}} \breve{J}_{1}(v_{1}(\cdot), u_{2}(\cdot)),\\
& \breve{J}_{2}(u_{1}(\cdot), u_{2}(\cdot)) = \min\limits_{v_{2}(\cdot)\in \mathscr{X}_{2}} \breve{J}_{2}(u_{1}(\cdot), v_{2}(\cdot)),
\end{aligned}
 \right.
\end{equation*}
subject to \eqref{mEq10}-\eqref{mEq11}.

 \textbf{Remark A.1}\,  Assume that two players partially
observe the state $(x, y^{v_{1}, v_{2}}, z_{1}^{v_{1}, v_{2}},
z_{2}^{v_{1}, v_{2}})$ and get the same observable information
$W_{2}$ in Example A.1. Similarly, we can formulate a non-zero sum
game of BSDE, and equivalently transform it into one with the same
Brownian motion observation, which is corresponding to the information structure (i) in Section \ref{PF}. The details of the deduction are omitted
for simplicity.

\textbf{Remark A.2}\, Recall the admissible control sets $\mathscr Y_1$ and
$\mathscr X_2$. Problem (A.1') is a non-zero sum stochastic
differential game of non-Markovian BSDE with asymmetric Brownian motion observation, which is corresponding to the information structure (ii) in Section \ref{PF}.

 \textbf{Example A.2.} Let the state and the cost
functional be same as \eqref{mEq8} and \eqref{mEq2}, respectively.
Suppose that $\left(x, y_{1}^{v_{1}, v_{2}}, z_{1}^{v_{1}, v_{2}},
z_{2}^{v_{1}, v_{2}}\right)$ is only partially observed by the
player $i$ $(i=1, 2)$ through
\begin{equation}\label{mEq3}
\left\{\begin{aligned}
  dW_{i}(t)=&\ \bar{h}_{i}\big(t, x(t)\big)dt+\sum_{j=1}^{2}\sigma_{ij}dw_{j}(t), \\
    W_{i}(0)=&\  0,
    \end{aligned}
    \right.
\end{equation}
respectively. Here $\bar{h}_{i}$ is uniformly bounded, and $\sigma=\begin{pmatrix}
                                                        \sigma_{11} & \sigma_{12} \\
                                                        \sigma_{21} & \sigma_{22} \\
                                                      \end{pmatrix}$ is an invertible constant matrix.
Admissible control set  for the player $i$ is defined by
$$\mathscr{V}_{i}=\left\{v_{i}(\cdot); v_{i}(t)\ \hbox{is}\ \mathscr{F}^{W_{i}}_{t}\hbox{-adapted and square integrable}\right\}.$$
Then the game problem is

  \textbf{Problem (A.2)}. Find a pair of admissible controls $(u_{1}, u_{2})$ such that
\begin{equation*}
 \left\{
 \begin{aligned}
& J_{1}(u_{1}(\cdot), u_{2}(\cdot)) = \min\limits_{v_{1}(\cdot)\in \mathscr{V}_{1}} J_{1}(v_{1}(\cdot), u_{2}(\cdot)),\\
& J_{2}(u_{1}(\cdot), u_{2}(\cdot)) = \min\limits_{v_{2}(\cdot)\in \mathscr{V}_{2}} J_{2}(u_{1}(\cdot), v_{2}(\cdot)),
\end{aligned}
 \right.
\end{equation*}
subject to \eqref{mEq8}, \eqref{mEq2} and \eqref{mEq3}.

 To simplify Problem (A.2), we set
\begin{align*}
&
\bar{h}(t, x(t))=
\begin{pmatrix}
  \bar{h}_{1}\big(t, x(t)\big) \\
 \bar{h}_{2}\big(t, x(t)\big)\\
\end{pmatrix},\quad  W(t)=\begin{pmatrix}
                                                             W_{1}(t) \\
                                                              W_{2}(t) \\
                                                            \end{pmatrix},\\
& \bar{\sigma}=\sigma^{-1}=
\begin{pmatrix}
  \bar{\sigma}_{11} & \bar{\sigma}_{12}\\
  \bar{\sigma}_{21} & \bar{\sigma}_{22} \\
\end{pmatrix},\\
&
\bar{c}(t, x(t))=
\begin{pmatrix}
  \bar{c}_{1}\big(t, x(t)\big) \\
  \bar{c}_{2}\big(t, x(t)\big)\\
\end{pmatrix}=\bar{\sigma} \bar{h}(t, x(t)),\\
 & \bar{w}(t)=w(t)+\int_{0}^{t}\bar{c}(s, x(s))ds,\\
 & \rho_{2}(t)=\exp\left\{-\int_{0}^{t}\bar{c}^{\,*}(s, x(s)) dw(s)-\frac{1}{2}\int_{0}^{t}|\bar{c}(s, x(s))|^{2}ds\right\}.
\end{align*}
Let $\frac{d\bar{\mathbb{P}}}{d\mathbb{P}}\Big |_{\mathscr{F}_{T}}=\rho_{2}(T).$ Similarly,
$\bar{\mathbb{P}}$ is a new probability measure, and consequently, $\bar{w}$ is a standard Brownian motion under $\bar{\mathbb{P}}$. Then
$$
dW(t)=\sigma d\bar{w}(t), \quad dw(t)=\bar{\sigma}dW(t)-\bar{c}(t, x(t))dt.
$$
We also set
$$X=x, Y^{v_{1}, v_{2}}=y^{v_{1}, v_{2}}, Z_{1}^{v_{1}, v_{2}}=\bar{\sigma}_{11}z_{1}^{v_{1}, v_{2}}+\bar{\sigma}_{21}z_{2}^{v_{1}, v_{2}},$$
$Z_{2}^{v_{1}, v_{2}}=\bar{\sigma}_{12}z_{1}^{v_{1}, v_{2}}+\bar{\sigma}_{22}z_{2}^{v_{1}, v_{2}}.$
With the notations above, \eqref{mEq8} is equivalently rewritten as
\begin{equation}\label{mEq4}
  \left\{
 \begin{aligned}
  -dY^{v_{1}, v_{2}}&(t)=\bar{g}\big(t, X(t), Y^{v_{1}, v_{2}}(t), Z_{1}^{v_{1}, v_{2}}(t), Z_{2}^{v_{1}, v_{2}}(t), v_{1}(t), v_{2}(t)\big)dt\\
            &\quad\quad -Z_{1}^{v_{1}, v_{2}}(t)dW_{1}(t)-Z_{2}^{v_{1}, v_{2}}(t)dW_{2}(t), \\
 Y^{v_{1}, v_{2}}&(T)=X(T)
 \end{aligned}
\right.
\end{equation}
with
$$
\label{mEq1}
  \left\{
 \begin{aligned}
 dX(t) = &\ \big(b\big(t, X(t)\big)-\delta_{1}(t)\bar{c}_{1}\big(t, X(t)\big)-\delta_{2}(t)\bar{c}_{2}\big(t, X(t)\big)\big)dt\\
 &\ +\left(\bar{\sigma}_{11}\delta_{1}(t)+\bar{\sigma}_{21}\delta_{2}(t)\right)dW_{1}(t)+\left(\bar{\sigma}_{12}\delta_{1}(t)+\bar{\sigma}_{22}\delta_{2}(t)\right)dW_{2}(t),\\
  X(0) = &\  0,
 \end{aligned}
\right.
$$
where
\begin{multline*}
  \bar{g}=g\left(t, Y^{v_{1}, v_{2}}(t), \frac{\bar{\sigma}_{22}Z_{1}^{v_{1}, v_{2}}(t)-\bar{\sigma}_{21}Z_{2}^{v_{1}, v_{2}}(t)}{\bar{\sigma}_{11}\bar{\sigma}_{22}-\bar{\sigma}_{12}\bar{\sigma}_{21}},\right.
 \left. \frac{\bar{\sigma}_{11}Z_{2}^{v_{1}, v_{2}}(t)-\bar{\sigma}_{12}Z_{1}^{v_{1}, v_{2}}(t)
}{\bar{\sigma}_{11}\bar{\sigma}_{22}-\bar{\sigma}_{12}\bar{\sigma}_{21}}, v_{1}(t), v_{2}(t)\right)\\
+\bar{c}_{1}(t, X(t))\frac{\bar{\sigma}_{22}Z_{1}^{v_{1}, v_{2}}(t)-\bar{\sigma}_{21}Z_{2}^{v_{1}, v_{2}}(t)
}{\bar{\sigma}_{11}\bar{\sigma}_{22}-\bar{\sigma}_{12}\bar{\sigma}_{21}}
+\bar{c}_{2}(t, X(t))\frac{\bar{\sigma}_{11}Z_{2}^{v_{1}, v_{2}}(t)-\bar{\sigma}_{12}Z_{1}^{v_{1}, v_{2}}(t)
}{\bar{\sigma}_{11}\bar{\sigma}_{22}-\bar{\sigma}_{12}\bar{\sigma}_{21}}.
\end{multline*}
Furthermore, we  assume that $\sigma$ is orthogonal in order to
guarantee that $W$ is also a standard Brownian motion under
$\bar{\mathbb{P}}$, under which \eqref{mEq4} is a non-Markov BSDE. On
the other hand, \eqref{mEq2} is rewritten as
\begin{multline}\label{mEq5}
  \bar{J}_{i}\big(v_{1}(\cdot), v_{2}(\cdot)\big)=\frac{1}{2} \mathbb{E}_{\bar{\mathbb{P}}}\left\{ \int_{0}^{T}\bar{l}_{i}\Big(t, X(t), Y^{v_{1}, v_{2}}(t), Z_{1}^{v_{1}, v_{2}}(t), Z_{2}^{v_{1}, v_{2}}(t), v_{1}(t), v_{2}(t)\Big)dt\right.\\
\left.+r_{i}\big(Y^{v_{1}, v_{2}}(0)\big)\vphantom{\int_{0}^{T}}\right\},
\end{multline}
where $\mathbb{E}_{\bar{\mathbb{P}}}$ denotes the expectation under $\bar{\mathbb{P}}$,
\begin{multline*}
  \bar{l}_{i}=\rho_{2}^{-1}(t)l_{i}\left(t,  Y^{v_{1}, v_{2}}(t), \frac{\bar{\sigma}_{22}Z_{1}^{v_{1}, v_{2}}(t)-\bar{\sigma}_{21}Z_{2}^{v_{1}, v_{2}}(t)}{\bar{\sigma}_{11}\bar{\sigma}_{22}-\bar{\sigma}_{12}\bar{\sigma}_{21}},\right.
 \left. \frac{\bar{\sigma}_{11}Z_{2}^{v_{1}, v_{2}}(t)-\bar{\sigma}_{12}Z_{1}^{v_{1}, v_{2}}(t)
}{\bar{\sigma}_{11}\bar{\sigma}_{22}-\bar{\sigma}_{12}\bar{\sigma}_{21}}, v_{1}, v_{2}\right)
\end{multline*}
and
$$
  \rho_{2}^{-1}(t)=\exp\left\{\int_{0}^{t}\bar{c}^{\,*}(s, X(s))\bar{\sigma}dW(s)-\frac{1}{2}\int_{0}^{t}|\bar{c}^{\,*}(s, X(s))|^{2}ds\right\}.
$$

Now Problem (A.2) is equivalently stated as follows.

\textbf{Problem (A.2')}. Find a pair of admissible controls $(u_{1}, u_{2})$ such that
\begin{equation*}
 \left\{
 \begin{aligned}
& \bar{J}_{1}(u_{1}(\cdot), u_{2}(\cdot)) = \min\limits_{v_{1}(\cdot)\in \mathscr{V}_{1}} \bar{J}_{1}(v_{1}(\cdot), u_{2}(\cdot)),\\
& \bar{J}_{2}(u_{1}(\cdot), u_{2}(\cdot)) = \min\limits_{v_{2}(\cdot)\in \mathscr{V}_{2}} \bar{J}_{2}(u_{1}(\cdot), v_{2}(\cdot)),
\end{aligned}
 \right.
\end{equation*}
subject to \eqref{mEq4}-\eqref{mEq5}.

 \textbf{Remark A.3}\, This is also a non-zero sum stochastic differential game of non-Markovian BSDE with mutually independent Brownian motion observation, which is corresponding to the information structure (iii) in Section \ref{PF}.






\end{document}